\documentclass{amsart}
\usepackage{amsfonts}

\usepackage{amscd}

\input{tcilatex}
\input tcilatex

\begin{document}
\author{Marek Kara\'{s}}
\title{There is no tame automorphism of $\Bbb{C}^{3}$ with multidegree $(3,4,5)$}
\date{}
\maketitle

\begin{abstract}
Let $F=(F_{1},\ldots ,F_{n}):\Bbb{C}^{n}\rightarrow \Bbb{C}^{n}$ be any
polynomial mapping. By multidegree of $F,$ denoted $\limfunc{mdeg}F,$ we
call the sequence of positive integers $(\deg F_{1},\ldots ,\deg F_{n}).$ In
this paper we addres the following problem: \textit{for which sequence }$%
(d_{1},\ldots ,d_{n})$\textit{\ there is an automorphism or tame
automorphism }$F:\Bbb{C}^{n}\rightarrow \Bbb{C}^{n}$\textit{\ with }$%
\limfunc{mdeg}F=(d_{1},\ldots ,d_{n}).$ We proved, among other things, that
there is no tame automorphism $F:\Bbb{C}^{3}\rightarrow \Bbb{C}^{3}$ with $%
\limfunc{mdeg}F=(3,4,5).$
\end{abstract}

\section{Introduction}

Let $F=(F_{1},F_{2}):\Bbb{C}^{2}\rightarrow \Bbb{C}^{2}$ be any polynomial
automorphism. By Jung van der Kullk theorem \cite{Jung,Kulk} we have that $%
\deg F_{1}|\deg F_{2}$ or $\deg F_{2}|\deg F_{1}.$ On the other hand if $%
d_{1},d_{2}$ are positive integers such that $d_{1}|d_{2}$ then $F=\Phi
_{2}\circ \Phi _{1},$ where 
\begin{eqnarray*}
\Phi _{1} &:&\Bbb{C}^{2}\ni (x,y)\mapsto (x+y^{d_{1}},y)\in \Bbb{C}^{2}, \\
\Phi _{2} &:&\Bbb{C}^{2}\ni (u,w)\mapsto (u,w+u^{\frac{d_{2}}{d_{1}}})\in 
\Bbb{C}^{2},
\end{eqnarray*}
is an automorphism of $\Bbb{C}^{2}$ such that $\limfunc{mdeg}%
F=(d_{1},d_{2}).\,$Similarily if $d_{2}|d_{1}$ we can writedown the
appropriate automorphism of $\Bbb{C}^{2}.$ Thus for the sequence of positive
integers $(d_{1},d_{2})$ to be the multidegree of some polynomial
automorphism of $\Bbb{C}^{2}$ is equivalent to satisy the condition: $%
d_{1}|d_{2}\,$or $d_{2}|d_{1}.$

It seems to be natural to ask for which sequence $(d_{1},\ldots ,d_{n})$
there is a polynomial automorphism $F:\Bbb{C}^{n}\rightarrow \Bbb{C}^{n}$
with $\limfunc{mdeg}F=(d_{1},\ldots ,d_{n}).$ Also, the question about
existanse of a tame automorphism $F:\Bbb{C}^{n}\rightarrow \Bbb{C}^{n}$ with 
$\limfunc{mdeg}F=(d_{1},\ldots ,d_{n})$ is natural. Recall that a tame
automorphism is, by definition, a composition of linear automorphisms and
triangular automorphisms, where a triangular automorphism is a mapping of
the following form 
\begin{equation*}
T:\Bbb{C}^{n}\ni \left\{ 
\begin{array}{l}
x_{1} \\ 
x_{2} \\ 
\vdots \\ 
x_{n}
\end{array}
\right\} \mapsto \left\{ 
\begin{array}{l}
x_{1} \\ 
x_{2}+f_{2}(x_{1}) \\ 
\vdots \\ 
x_{n}+f_{n}(x_{1},\ldots ,x_{n-1})
\end{array}
\right\} \in \Bbb{C}^{n}.
\end{equation*}
By $\limfunc{Tame}(\Bbb{C}^{n})$ we will denote the group of all tame
automorphimsm of $\Bbb{C}^{n}.$ This is, of course, a subgroup of the group $%
\limfunc{Aut}(\Bbb{C}^{n})$ of all polynomial automorphisms of $\Bbb{C}^{n}.$

It is easy to see that if there is an automorphism (or tame automorphism) $F:%
\Bbb{C}^{n}\rightarrow \Bbb{C}^{n}$ such that $\limfunc{mdeg}F=(d_{1},\ldots
,d_{n})$ then there is, also, an automorphism (or tame automorphism) $%
\widetilde{F}:\Bbb{C}^{n}\rightarrow \Bbb{C}^{n}$ such that $\limfunc{mdeg}%
\widetilde{F}=(d_{\sigma (1)},\ldots ,d_{\sigma (n)})$ for any permutation $%
\sigma \,$of the set $\{1,\ldots ,n\}.$ Thus in our considerations, without
loose of generality, we can assume that $d_{1}\leq d_{2}\leq \ldots \leq
d_{n}.$

\section{Some simple remarks}

In this section we make some simple but useful remarks about existense of
automorphism and tame automorphism with given multidegree.

\begin{proposition}
If for $1\leq d_{1}\leq \ldots d_{n}$ there is a sequence of integers $1\leq
i_{1}<\ldots <i_{m}\leq n,$ with $m<n,$ such that there exists an
automorphism $G$ of $\Bbb{C}^{m}$ with $\limfunc{mdeg}G=(d_{i_{1}},\ldots
,d_{i_{m}}),$ then there exists an automorphis $F$ of $\Bbb{C}^{n}$ with $%
\limfunc{mdeg}F=(d_{1},\ldots ,d_{n}).$ Moreover, if we assume that $G$ is a
tame automorphism, then there is a tame automorphism $F$ of $\Bbb{C}^{n}$
such that $\limfunc{mdeg}F=(d_{1},\ldots ,d_{n}).$
\end{proposition}

\begin{proof}
Let $j_{1},\ldots ,j_{n-m}\in \Bbb{N}$ be such that $1\leq j_{1}<\ldots
<j_{n-m}\leq n$ and $\{i_{1},\ldots ,i_{m}\}\cup \{j_{1},\ldots
,j_{n-m}\}=\{1,\ldots ,n\}.$ In this situation we have, of course, $%
\{i_{1},\ldots ,i_{m}\}\cap \{j_{1},\ldots ,j_{n-m}\}=\emptyset .$ Consider
the mapping $h=(h_{1},\ldots ,h_{n}):\Bbb{C}^{n}\rightarrow \Bbb{C}^{n}$
given by the formulas 
\begin{equation*}
h_{k}(x_{1},\ldots ,x_{n})=\left\{ 
\begin{array}{ll}
x_{k} & \text{for }k\in \{i_{1},\ldots ,i_{m}\} \\ 
x_{k}+(x_{i_{1}})^{d_{k}}\quad & \text{for }k\in \{j_{1},\ldots ,j_{n-m}\}
\end{array}
\right. .
\end{equation*}
Of course $h$ is an automorphism of $\Bbb{C}^{n}$ and $\deg h_{k}=d_{k}$ for 
$k\in \{i_{1},\ldots ,i_{m}\}.$

Consider, also, the mapping $g=(g_{1},\ldots ,g_{n}):\Bbb{C}^{n}\rightarrow 
\Bbb{C}^{n}$ given by the formulas 
\begin{equation*}
g_{k}(u_{1},\ldots ,u_{n})=\left\{ 
\begin{array}{ll}
G_{l}(u_{i_{1}},\ldots ,u_{i_{m}})\quad & \text{for }k=i_{l} \\ 
u_{k}\quad & \text{for }k\in \{j_{1},\ldots ,j_{n-m}\}
\end{array}
\right. .
\end{equation*}
It is easy to see that $g$ is an automorphism of $\Bbb{C}^{n}$ and $\deg
g_{k}=d_{k}$ for $k\in \{j_{1},\ldots ,j_{n-m}\}.$

Now taking $F=g\circ h$ we obtain an automorphism of $\Bbb{C}^{n}$ such that 
$\deg F_{i}=d_{i}$ for all $i\in \{1,\ldots ,n\}.$
\end{proof}

\begin{proposition}
\label{prop_sum_d_i}If for a sequence of integers $1\leq d_{1}\leq \ldots
\leq d_{n}$ there is $i\in \{1,\ldots ,n\}$ such that 
\begin{equation*}
d_{i}=\sum_{j=1}^{i-1}k_{j}d_{j}\qquad \text{with }k_{j}\in \Bbb{N},
\end{equation*}
then there exists a tame automorphism $F$ of $\Bbb{C}^{n}$ with $\limfunc{%
mdeg}F=(d_{1},\ldots ,d_{n}).$
\end{proposition}

\begin{proof}
Consider the following two mappings $h=(h_{1},\ldots ,h_{n}):\Bbb{C}%
^{n}\rightarrow \Bbb{C}^{n}$ and $g=(g_{1},\ldots ,g_{n}):\Bbb{C}%
^{n}\rightarrow \Bbb{C}^{n}$ given by the formulas 
\begin{equation*}
h_{k}(x_{1},\ldots ,x_{n})=\left\{ 
\begin{array}{ll}
x_{k}\quad & \text{for }k=i \\ 
x_{k}+x_{i}^{d_{k}}\quad & \text{for }k\neq i
\end{array}
\right.
\end{equation*}
and 
\begin{equation*}
g_{k}(u_{1},\ldots ,u_{n})=\left\{ 
\begin{array}{ll}
u_{k}+u_{1}^{k_{1}}\cdots u_{i-1}^{k_{i-1}}\quad & \text{for }k=i \\ 
u_{k}\quad & \text{for }k\neq i
\end{array}
\right. .
\end{equation*}
Now it is easy to see that $h$ and $g$ are automorphisms of $\Bbb{C}^{n}$
such that $\deg h_{k}=d_{k}$ for $k\neq i$ and $\deg g_{i}=d_{i}.$ Since,
also, $h_{i}(x_{1},\ldots ,x_{n})=x_{i}$ and $g_{k}(u_{1},\ldots
,u_{n})=u_{k}$ for $k\neq i,$ then it is easy to see that for the
automorphism $F=g\circ h$ we have $\deg F_{k}=d_{k}$ for all $k\in
\{1,\ldots ,n\}.$
\end{proof}

\begin{corollary}
\label{cor_male_d1}If for a sequence of integers $1\leq d_{1}\leq \ldots
\leq d_{n}$ we have $d_{1}\leq n-1,$ then there exists a tame automorphis $F$
of $\Bbb{C}^{n}$ with $\limfunc{mdeg}F=(d_{1},\ldots ,d_{n}).$
\end{corollary}

\begin{proof}
Let $r_{i}\in \{0,1,\ldots ,d_{1}-1\},$ for $i=2,\ldots ,n,$ be such that $%
d_{i}\equiv r_{i}(\func{mod}d_{1}),$ for $i=2,\ldots ,n.$ If there is an $%
i\in \{2,\ldots ,n\}$ such that $r_{i}=0,$ then $d_{i}=kd_{1}$ for some $%
k\in \Bbb{N}\backslash \{0\}$ and by Proposition \ref{prop_sum_d_i}, there
exists an automorphis $F$ of $\Bbb{C}^{n}$ with the desired properties.

Thus we can assume that $r_{i}\neq 0$ for all $i=2,\ldots ,n.$ Since $%
d_{1}-1<n-1,$ then there are $i,j\in \{2,\ldots ,n\},i\neq j,$ such that $%
r_{i}=r_{j}.$ Without lose of generality we can assume that $i<j.$ In this
situation we have $d_{j}=d_{i}+kd_{1}$ for some $k\in \Bbb{N}.$ Then by
Proposition \ref{prop_sum_d_i} there exists an automorphis $F$ of $\Bbb{C}%
^{n}$ with the desired properies.
\end{proof}

\section{Examples\label{sec_exmples}}

In this section we give some positive results about existence of tame
automorphisms of $\Bbb{C}^{3}$ with given multidegree $(d_{1},d_{2},d_{3}).$
The first one is the following.

\begin{example}
For every $d_{2},d_{3}\in \Bbb{N},2\leq d_{2}\leq d_{3},$ there is a tame
automorphis $F$ of $\Bbb{C}^{3}$ such that 
\begin{equation*}
\limfunc{mdeg}F=(2,d_{2},d_{3}).
\end{equation*}
This is a consequence of Corollary \ref{cor_male_d1}.
\end{example}

\begin{example}
For any $d_{3}\geq 4$ such that $d_{3}\neq 5$ there is a tame automorphis $F$
of $\Bbb{C}^{3}$ such that 
\begin{equation*}
\limfunc{mdeg}F=(3,4,d_{3}).
\end{equation*}
\end{example}

\begin{proof}
We have 
\begin{equation*}
4=0\cdot 3+1\cdot 4
\end{equation*}
and 
\begin{equation*}
d_{3}=\left\{ 
\begin{array}{ll}
(2+k)\cdot 3+0\cdot 4 & \text{for }d_{3}=6+3k \\ 
(1+k)\cdot 3+1\cdot 4 & \text{for }d_{3}=7+3k \\ 
(0+k)\cdot 3+2\cdot 4 & \text{for }d_{3}=8+3k
\end{array}
\right. .
\end{equation*}
Thus we can apply Proposition \ref{prop_sum_d_i}.
\end{proof}

\begin{example}
For any $d_{3}\geq 5$ such that $d_{3}\neq 7$ there is a tame automorphis $F$
of $\Bbb{C}^{3}$ such that 
\begin{equation*}
\limfunc{mdeg}F=(3,5,d_{3}).
\end{equation*}
\end{example}

\begin{proof}
We have 
\begin{equation*}
5=0\cdot 3+1\cdot 5,\qquad 6=2\cdot 3+0\cdot 5
\end{equation*}
and 
\begin{equation*}
d_{3}=\left\{ 
\begin{array}{ll}
(1+k)\cdot 3+1\cdot 5 & \text{for }d_{3}=8+3k \\ 
(3+k)\cdot 3+0\cdot 5 & \text{for }d_{3}=9+3k \\ 
(0+k)\cdot 3+2\cdot 5 & \text{for }d_{3}=10+3k
\end{array}
\right. .
\end{equation*}
Thus we can apply Proposition \ref{prop_sum_d_i}.
\end{proof}

\begin{example}
For any $d_{3}\geq 5$ such that $d_{3}\neq 6,7,11$ there is a tame
automorphism $F$ of $\Bbb{C}^{3}$ such that 
\begin{equation*}
\limfunc{mdeg}F=(4,5,d_{3}).
\end{equation*}
\end{example}

\begin{proof}
We have 
\begin{eqnarray*}
5 &=&0\cdot 4+1\cdot 5,\qquad 8=2\cdot 4+0\cdot 5, \\
9 &=&1\cdot 4+1\cdot 5,\qquad 10=0\cdot 4+2\cdot 5
\end{eqnarray*}
and 
\begin{equation*}
d_{3}=\left\{ 
\begin{array}{ll}
(3+k)\cdot 4+0\cdot 5 & \text{for }d_{3}=12+4k \\ 
(2+k)\cdot 4+1\cdot 5 & \text{for }d_{3}=13+4k \\ 
(1+k)\cdot 4+2\cdot 5 & \text{for }d_{3}=14+4k \\ 
(0+k)\cdot 4+3\cdot 5 & \text{for }d_{3}=15+4k
\end{array}
\right. .
\end{equation*}
Thus we can apply Proposition \ref{prop_sum_d_i}.
\end{proof}

The above examples justifies the following question.

\begin{description}
\item[Quastion]  Is there any automorphism (or tame automorphism) $F$ of $%
\Bbb{C}^{3}$ such that 
\begin{equation*}
\limfunc{mdeg}F\in \{(3,4,5),(3,5,7),(4,5,6),(4,5,7),(4,5,11)\}?
\end{equation*}
\end{description}

\section{Partial answer}

In this section we give partial answer for the quastion established in the
last section. Namely we show the following

\begin{theorem}
\label{tw_345}There is no tame automorphism $F=(F_{1},F_{2},F_{3})$ of $\Bbb{%
C}^{3}$ such that 
\begin{equation*}
\limfunc{mdeg}F=(3,4,5).
\end{equation*}
\end{theorem}

Before we make a proof of Theorem \ref{tw_345} we recall some results and
notions from the papers of Shestakov and Umirbayev \cite{sh umb1,sh umb2}.

\begin{definition}
\textit{(\cite{sh umb1}, Definition 1) }\label{def_*-red}A pair $f,g\in
k[X_{1},\ldots ,X_{n}]$ is called *-reduced if\newline
(i) $f,g$ are algebraically independnt;\newline
(ii) $\overline{f},\overline{g}$ are algebraically dependnt, where $%
\overline{h}$ denotes the highest homogeneous part of $h$;\newline
(iii) $\overline{f}\notin [\overline{g}]$ and $\overline{g}\notin k[%
\overline{f}].$
\end{definition}

\begin{definition}
\textit{(\cite{sh umb1}, Definition 1) }Let $f,g\in k[X_{1},\ldots ,X_{n}]$
be a *-reduced pair with $\deg f<\deg g.$ Put $p=\frac{\deg f}{\gcd (\deg
f,\deg g)}.$ In this sitation the pair $f,g$ is called $p-$reduced pair.
\end{definition}

\begin{theorem}
\textit{(\cite{sh umb1}, Theorem 2)}\label{tw_deg_g_fg} Let $f,g\in
k[X_{1},\ldots ,X_{n}]$ be a $p-$reduced pair, and let $G(x,y)\in k[x,y]$
with $\deg _{y}G(x,y)=pq+r,0\leq r<p.$ Then 
\begin{equation*}
\deg G(f,g)\geq q\left( p\deg g-\deg g-\deg f+\deg [f,g]\right) +r\deg g.
\end{equation*}
\end{theorem}

In the above theorem $[f,g]$ means the Poisson bracket of $f$ and $g,$ but
for us it is only important that 
\begin{equation*}
\deg [f,g]=2+\underset{1\leq i<j\leq n}{\max }\deg \left( \frac{\partial f}{%
\partial x_{i}}\frac{\partial g}{\partial x_{j}}-\frac{\partial f}{\partial
x_{j}}\frac{\partial g}{\partial x_{i}}\right)
\end{equation*}
if $f,g$ are algebraically independent, and $\deg [f,g]=0$ if $f,g$ are
algebraically dependent.

Notice, also, that the estimation from Theorem \ref{tw_deg_g_fg} is true
even if the condition (ii) of Definition \ref{def_*-red} is not satisfied.
Indeed, if $G(x,y)=\sum_{i,j}a_{i,j}x^{i}y^{j},$ then, by the algebraic
independence of $\overline{f}$ and $\overline{g}$ we have: 
\begin{eqnarray*}
\deg G(f,g) &=&\underset{i,j}{\max }\deg (a_{i,j}f^{i}g^{j})\geq \deg
_{y}G(x,y)\cdot \deg g= \\
&=&(qp+r)\deg g\geq q(p\deg g-\deg f-\deg g+\deg [f,g])+r\deg g.
\end{eqnarray*}
The last inequality is a consequence of the fact that $\deg [f,g]\leq \deg
f+\deg g.$

We will also use the following theorem.

\begin{theorem}
\label{tw_type_1-4}\textit{(\cite{sh umb1}, Theorem 3) }Let $%
F=(F_{1},F_{2},F_{3})\,$be a tame automorphism of $\Bbb{C}^{3}.$ If $\deg
F_{1}+\deg F_{2}+\deg F_{3}>3$ (in other words if $F$ is not a lienar
automorphism), then $F$ admits either an elementary reduction or a reduction
of types I-IV (see \cite{sh umb1} Definitions 2-4).
\end{theorem}

Let us recall that an automorphism $F=(F_{1},F_{2},F_{3})$ admits an
elementary reduction if there exists a polynomial $g\in \Bbb{C}[x,y]$ and a
permutation $\sigma $ of the set $\{1,2,3\}$ such that $\deg (F_{\sigma
(1)}-g(F_{\sigma (2)},F_{\sigma (3)}))<\deg F_{\sigma (1)}.$

Now we are in a position to prove Theorem \ref{tw_345}

\begin{proof}
\textit{(of Theorem \ref{tw_345}) }Assume that $F=(F_{1},F_{2},F_{3})$ is an
automorphism of $\Bbb{C}^{3}$ such that $\limfunc{mdeg}F=(3,4,5).$ We will
show that this hypotethical automorphism (we do not know if there is any)
can not be tame. First of all, notice that any pair $F_{i},F_{j}$ with $%
i,j\in \{1,2,3\},i\neq j,$ satisfies the conditions (i) and (iii) of
Definition \ref{def_*-red}. Indeed, it follows by the fact that $%
F_{1},F_{2},F_{3}$ are algebraically independent and that $3,4\notin 5\Bbb{N}%
,3,5\notin 4\Bbb{N}$ and $4,5\notin 3\Bbb{N}.$ By Theorem \ref{tw_type_1-4}
it is enough to show that $F$ does not admit neither reductions of type I-IV
nor elementary reduction.

By a contrary, assume that $(F_{1},F_{2},F_{3})$ admits a reduction of type
I or II. Then by the definition (see \cite{sh umb1} Definition 2 and 3), for
some number $n\in \Bbb{N}\backslash \{0\}$ and some permutation $\sigma $ of
the set $\{1,2,3\}$ we have $\deg F_{\sigma (1)}=2n$ and $\deg F_{\sigma
(2)}=ns,$ where $s\geq 3$ is an odd number. But in the sequence $3,4,5$
there is only one even number, namely $4.$ Thus $2n=4,n=2$ and then $ns$ is,
also, an even number, a contradiction.

Now assume, by a contrary, that $(F_{1},F_{2},F_{3})$ admits a reduction of
type III or IV. Then by the definition (see \cite{sh umb1} Definition 4),
for some number $n\in \Bbb{N}\backslash \{0\}$ and some permutation $\sigma $
of the set $\{1,2,3\}$ we have $\deg F_{\sigma (1)}=2n$ and either 
\begin{equation}
\deg F_{\sigma (2)}=3n,n<\deg F_{\sigma (3)}\leq \tfrac{3}{2}n  \label{row_1}
\end{equation}
or 
\begin{equation}
\tfrac{5}{2}n<\deg F_{\sigma (2)}\leq 3n,\deg F_{\sigma (3)}=\tfrac{3}{2}n.
\label{row_2}
\end{equation}
Of course, as before, we have $2n=4,n=2.$ Since $3n=6,$ then (\ref{row_1})
is impossible, and since $\tfrac{5}{2}n=5,3n=6$ and $\deg F_{\sigma (2)}\in 
\Bbb{N},$ then (\ref{row_2}) is impossible. Thus we obtain a contradiction.

Thus, in order to show that $(F_{1},F_{2},F_{3})$ can not be a tame
automorphism, by Theorem \ref{tw_type_1-4}, it is enough to show that $%
(F_{1},F_{2},F_{3})$ does not admit an elemnetary reduction.

Let us assume that 
\begin{equation*}
(F_{1},F_{2},F_{3}-g(F_{1},F_{2})),
\end{equation*}
where $g\in k[x,y],$ is an elementary reduction of $(F_{1},F_{2},F_{3}).$
Thus, in particular, we have $\deg g(F_{1},F_{2})=5.$ But it is impossible.
Indded, by Theorem \ref{tw_deg_g_fg}, we have 
\begin{equation}
\deg g(F_{1},F_{2})\geq q(pm-m-n+\deg [F_{1},F_{2}])+mr,  \label{row1}
\end{equation}
where $n=\deg F_{1},m=\deg F_{2},p=n/\limfunc{GCD}(n,m)$ and $\deg
_{y}g(x,y)=qp+r$ with $0\leq r<p.\,$In our case we have $n=3,m=4,p=3.$ Since 
$F_{1},F_{2}$ are algebraically independent, $\deg [F_{1},F_{2}]\geq 2.$
Thus (\ref{row1}) can be rewriten as follows 
\begin{equation*}
\deg g(F_{1},F_{2})\geq q(3\cdot 4-4-3+\deg [F_{1},F_{2}])+4r.
\end{equation*}
Since, also, $3\cdot 4-4-3+\deg [F_{1},F_{2}]=5+$ $\deg [F_{1},F_{2}]\geq
7>5,$ then $q$ must be zero, and $r$ must be not greater than $1.$ This
means that $g(F_{1},F_{2})=g_{1}(F_{1})+g_{2}(F_{1})F_{2}\,$for some $%
g_{1},g_{2}\in k[x].$ Since $3\Bbb{N\cap (}4\Bbb{+}3\Bbb{\Bbb{N})=\emptyset }%
,$ then $\deg g(F_{1},F_{2})=\max \{3\deg g_{1},4+3\deg g_{2}\}.$ But, since 
$5\notin 3\Bbb{N}\cup (4+3\Bbb{N}),$ then we obtain a contradiction.

Now, let us assume that 
\begin{equation*}
(F_{1},F_{2}-g(F_{1},F_{3}),F_{3}),
\end{equation*}
where $g\in k[x,y],$ is an elementary reduction of $(F_{1},F_{2},F_{3}).$ In
this case we have $\deg g(F_{1},F_{3})=4.$ But it is impossible. Indeed, by
Theorem \ref{tw_deg_g_fg}, we have 
\begin{equation*}
\deg g(F_{1},F_{3})\geq q(pm-m-n+\deg [F_{1},F_{2}])+mr,
\end{equation*}
where $n=3,m=5,p=3$ and $\deg _{y}g(x,y)=3q+r$ with $0\leq r<3.$ Since $%
pm-m-n+\deg [F_{1},F_{2}]=7+\deg [F_{1},F_{2}]>4,$ then $q$ must be zero.
Also, $r$ must be zero, because $m=5>4.$ Thus $g(F_{1},F_{3})=g(F_{1}),$ and
then $\deg g(F_{1},F_{3})=3\deg g.$ Since $4\notin 3\Bbb{N},$ then we obtain
a contradiction.

And finally, let us assume that 
\begin{equation*}
(F_{1}-g(F_{2},F_{3}),F_{2},F_{3}),
\end{equation*}
where $g\in k[x,y],$ is an elementary reduction of $(F_{1},F_{2},F_{3}).$
Similarilly, as before, we obtain 
\begin{equation*}
3=\deg g(F_{2},F_{3})\geq q(4\cdot 5-5-4+\deg [F_{2},F_{3}])+5r,
\end{equation*}
where $\deg _{y}g(x,y)=4q+r$ with $0\leq r<4.$ Then $q$ and $r$ must be
zero. Thus $g(F_{2},F_{3})=g(F_{2}),$ and then $\deg g(F_{2},F_{3})=4\deg g.$
Since $3\notin 4\Bbb{N},$ then we obtain a contradiction.
\end{proof}

In the similar way we can show the following theorem.

\begin{theorem}
There is no tame automorphism $F$ of $\Bbb{C}^{3}$ such that 
\begin{equation*}
\limfunc{mdeg}F\in \{(3,5,7),(4,5,7),(4,5,11)\}.
\end{equation*}
\end{theorem}

By above theorem, Theorem \ref{tw_345}, Corollary \ref{cor_male_d1} and
examples from section \ref{sec_exmples} we have the following theorem.

\begin{theorem}
In the following statements $\limfunc{mdeg}$ is considered as a map from the
set of all endomorphisms of $\Bbb{C}^{n}$ into the set $\Bbb{N}^{n}.$\newline
(i) For all integers $d_{3}\geq d_{2}\geq 2,$ $(2,d_{2},d_{n})\in \limfunc{%
mdeg}(\limfunc{Tame}(\Bbb{C}^{3})).$\newline
(ii) If $d_{3}\geq 4,$ then $(3,4,d_{3})\in \limfunc{mdeg}(\limfunc{Tame}(%
\Bbb{C}^{3}))$ if and only if $d_{3}\neq 5.$\newline
(iii) If $d_{3}\geq 5,$ then $(3,5,d_{3})\in \limfunc{mdeg}(\limfunc{Tame}(%
\Bbb{C}^{3}))$ if and only if $d_{3}\neq 7.$\newline
(iv) If $d_{3}\geq 5$ and $d_{3}\neq 6,7,11,$ then $(4,5,d_{3})\in \limfunc{%
mdeg}(\limfunc{Tame}(\Bbb{C}^{3})).$\newline
\end{theorem}

\vspace{1cm}

\textsc{Marek Kara\'{s}\newline
Instytut Matematyki\newline
Uniwersytetu Jagiello\'{n}skiego\newline
ul. \L ojasiewicza 6}\newline
\textsc{30-348 Krak\'{o}w\newline
Poland\newline
} e-mail: Marek.Karas@im.uj.edu.pl

\end{document}